\documentclass{au}
\usepackage{amssymb, natbib, latexcad}

\theoremstyle{definition}
\newtheorem{definition}{Definition}[section]
\newtheorem{example}[definition]{Example}

\theoremstyle{plain}

\newtheorem{lemma}[definition]{Lemma}

\newtheorem{proposition}[definition]{Proposition}

\theoremstyle{remark}

\newtheorem{remark}[definition]{Remark}

\newcommand{\jumpname}{\Gamma}              
\newcommand{\aut}[1]{\mathrm{Aut}(#1)}      
\newcommand{\mbf}[1]{\mathbf{#1}}           
\newcommand{\spn}[1]{\langle#1\rangle}      
\newcommand{\probname}{p}                   
\newcommand{\prob}[2]{\probname_{#1}(#2)}   
\newcommand{\genname}{\mathrm{Gen}}         
\newcommand{\gen}[2]{\genname_{#1}(#2)}     
\def\jump[#1]#2#3{\jumpname_{#1}(#2,\,#3)}  
\def\orbit[#1]{O_{#1}}                      
\def\orbeq[#1]{\sim_{#1}}                   
\def\orbclass[#1]#2{#2_{\orbeq[#1]}}        
\newcommand{\neutral}{e}                    
\newcommand{\cyclic}[1]{C_{#1}}             
\newcommand{\hassename}{\mathcal H}         
\newcommand{\glbidx}[2]{\hassename_{#2}(#1)}
\newcommand{\isoidx}[3]{
    \hassename_{#3}(#1|#2)}                 
\newcommand{\orbidx}[3]{
    \hassename_{#3}^*(#1|#2)}               
\newcommand{\paige}[1]{M^*(#1)}             
\newcommand{\order}[1]{|#1|}                

\title[Random generators of given orders]
{Random Generators of Given Orders and the Smallest Simple Moufang Loop}
\author{Petr Vojt\v echovsk\'y}
\thanks{\\ While working on this paper the author has been partially
supported by Grant Agency of Charles University, grant number
269/2001/B-MAT/MFF}
\address{Department of Mathematics, Iowa State University, Ames, IA 50011, U.S.A.}
\email{petr@math.du.edu}



\begin{document}


\begin{abstract}
The probability that $m$ randomly chosen elements of a finite power associative
loop $C$ have prescribed orders and generate $C$ is calculated in terms of
certain constants $\jumpname$ related to the action of $\aut{C}$ on the subloop
lattice of $C$. As an illustration, all meaningful probabilities of random
generation by elements of given orders are found for the smallest
nonassociative simple Moufang loop.
\end{abstract}

\maketitle

\section{Random generators of given orders}\label{Sc:Random}

\noindent Let $C$ be a power associative loop, i.e., a loop where every element
generates a group. Given an $m$-tuple $\mbf{a}=(a_0$, $\dots$, $a_{m-1})$ of
elements of $C$, let $\mbf{A}=\{A_i\}_{i=0}^{m}$ be the sequence of nested
subloops $A_i\le C$ such that $A_0$ is the smallest subloop of $C$, and
$A_{i+1}=\spn{A_i,\,a_i}$. Note that $A_m$ is independent of the order of the
elements $a_0$, $\dots$, $a_{m-1}$ in $\mbf{a}$.

Denote by $\gen{m}{C}$ the set of all $m$-tuples $\mbf{a}\in C^m$ with $A_m=C$.
Then the probability that $m$ randomly chosen elements of $C$ generate $C$ is
\begin{equation}\label{Eq:Prob}
    \prob{m}{C}=|C|^{-m}\cdot|\gen{m}{C}|.
\end{equation}

This notion can be refined in a natural way. For $1\le i\le n$, let $D_i$ be
the set of all elements of $C$ of order $i$. Two $m$-tuples of integers
$\mbf{r}=(r_0$, $\dots$, $r_{m-1})$, $\mbf{s}=(s_0$, $\dots$, $s_{m-1})$ are
said to be of \emph{the same type} if $r_0$, $\dots$, $r_{m-1}$ is a
permutation of $s_0$, $\dots$, $s_{m-1}$. We say that $\mbf{a}=(a_0$, $\dots$,
$a_{m -1})\in C^m$ is of \emph{type} $\mbf{r}$ if there is $\mbf{s}=(s_0$,
$\dots$, $s_{m-1})$ of the same type as $\mbf{r}$ satisfying $a_i\in D_{s_i}$,
for $0\le i\le m-1$.

Let $\gen{\mbf{r}}{C}\subseteq \gen{m}{C}$ be the set of generating $m$-tuples
of type $\mbf{r}$. Then
\begin{equation}\label{Eq:TProb}
    \prob{\mbf{r}}{C}=|C|^{-m}\cdot |\gen{\mbf{r}}{C}|
\end{equation}
is the probability that $m$ randomly chosen elements $a_0$, $\dots$,
$a_{m-1}\in C$ generate $C$ and $(a_0$, $\dots$, $a_{m-1})$ is of type
$\mbf{r}$.

For $A$, $B\le C$ and an integer $i$, let $\jump[i]{A}{B}$ be the cardinality
of the set of all elements $x\in D_i$ such that $\spn{A,\,x}\in\orbit[B]$,
where $\orbit[B]$ is the orbit of $B$ under the natural action of $\aut{C}$ on
the subloop lattice of $C$. Also, let $\jump[]{A}{B}$ be the cardinality of the
set of all elements $x\in C$ such that $\spn{A,\,x}\in\orbit[B]$.

We are going to divide $C^m$ into certain equivalence classes. Two $m$-tuples
$\mbf{a}$, $\mbf{b}\in C^m$ with associated nested subloops $\{A_i\}_{i=0}^m$,
$\{B_i\}_{i=0}^m$ will be called \emph{orbit-equivalent} if
$A_i\in\orbit[B_i]$, for $0\le i\le m$. We write $\mbf{a}\orbeq[]\mbf{b}$.

The size of the equivalence class $\orbclass[]{\mbf{a}}$ is easy to calculate
with help of the constants $\jump[]{A}{B}$. There are $\jump[]{A_0}{A_1}$
elements $x$ such that $\spn{A_0,\,x}\in\orbit[A_1]$. Once we are in the orbit
$\orbit[A_i]$, we can continue on the way to $\orbit[A_{i+1}]$ by adding one of
the $\jump[]{A_i}{A_{i+1}}$ elements $x_{i+1}$ to $\spn{x_0,\,\dots,\,x_i}$.
Thus,
\begin{equation}\label{Eq:Orbclass}
    |\orbclass[]{\mbf{a}}|=\prod_{i=0}^{m-1}\jump[]{A_i}{A_{i+1}}.
\end{equation}
Since
\begin{equation}\label{Eq:Gen}
    |\gen{m}{C}|=\sum_{\orbclass[]{\mbf{a}}\in\gen{m}{C}/\orbeq[]}
        |\orbclass[]{\mbf{a}}|,
\end{equation}
we can combine $(\ref{Eq:Prob})$, $(\ref{Eq:Orbclass})$ and $(\ref{Eq:Gen})$ to
obtain
\begin{equation}\label{Eq:ExplicitProb}
    \prob{m}{C}=|C|^{-m}\sum_{\orbclass[]{\mbf{a}} \in \gen{m}{C}/\orbeq[]}
        \;\prod_{i=0}^{m-1}\jump[]{A_i}{A_{i+1}}.
\end{equation}

\begin{example}
Let us illustrate $(\ref{Eq:ExplicitProb})$ by calculating the probability that
two randomly chosen elements of $S_3$ generate the entire symmetric group. Let
$\neutral$ be the neutral element of $S_3$. There are $3$ subgroups isomorphic
to $\cyclic{2}$ (all in one orbit of transitivity) and a unique subgroup
isomorphic to $\cyclic{3}$. As $\jump[]{\neutral}{\cyclic{2}}=3$,
$\jump[]{\neutral}{\cyclic{3}}=2$, $\jump[]{\cyclic{2}}{S_3}=4$ and
$\jump[]{\cyclic{3}}{S_3}=3$, we have $\prob{2}{S_3}=(3\cdot 4+2\cdot
3)/36=1/2$, as expected.
\end{example}

Write $\orbeq[\mbf{r}]$ for the restriction of $\orbeq[]$ onto
$\gen{\mbf{r}}{C}$, and observe that
\begin{equation}\label{Eq:Torbclass}
    |\orbclass[\mbf{r}]{\mbf{a}}|=\sum_{\mbf{s}=(s_0,\,\dots,\,s_{m-1})}
        \;\prod_{i=0}^{m-1}\jump[s_i]{A_i}{A_{i+1}},
\end{equation}
where the summation runs over all $m$-tuples $\mbf{s}$ of the same type as
$\mbf{r}$. Consequently,
\begin{equation}\label{Eq:Tprob}
    \prob{\mbf{r}}{C}=|C|^{-m}
        \sum_{\orbclass[\mbf{r}]{\mbf{a}}\in\gen{\mbf{r}}{C}/\orbeq[\mbf{r}]}
        \;\;\sum_{\mbf{s}=(s_0,\,\dots,\,s_{m-1})}
        \;\prod_{i=0}^{m-1}\jump[s_i]{A_i}{A_{i+1}}.
\end{equation}

\begin{remark}
All concepts of this section can be generalized to any finite universal algebra
$C$ with subsets $D_i$ closed under the action of $\aut{C}$.
\end{remark}

\section{Random generators of given orders for $M^*(2)$}

\noindent We assume from now on that the reader is familiar with the notation
and terminology of \cite{VojtechovskyHasse}.

The value of $\jump[i]{A}{B}$ in (\ref{Eq:Tprob}) can be calculated with help
of Hasse constants, provided $A$ is maximal in $B$. Namely,
\begin{equation}\label{Eq:Hasse}
    \jump[i]{A}{B}=\orbidx{A}{B}{C}\cdot|D_i\cap(B\setminus A)|.
\end{equation}
This is obvious because $\orbidx{A}{B}{C}$ counts the number of subloops of $C$
containing $A$ and in the same orbit as $B$.

By $(\ref{Eq:Tprob})$ and $(\ref{Eq:Hasse})$, we should be able to calculate
$\prob{m}{C}$, $\prob{\mbf{r}}{C}$ when all Hasse constants for $C$ are known
and when the lattice of subloops of $C$ is not too high. These probabilities
can be used alongside order statistics to recognize black box loops, for
instance (cf.\ \cite{GroupsAndComputationIII}).

Building on the results of \cite{VojtechovskyHasse} substantially, we proceed
to calculate all meaningful probabilities $\prob{m}{C}$, $\prob{\mbf{r}}{C}$
for $C=\paige{2}$---the smallest nonassociative simple Moufang loop.

All Hasse constants for $C$ are summarized in \cite[Fig.\
$1$]{VojtechovskyHasse}, so we can easily evaluate all constants
$\jump[i]{A}{B}$ with $A$ maximal in $B$. For example, since
$\orbidx{E_8}{M(A_4)}{C}=3$ and $|D_2\cap(M(A_4)\setminus E_8)|=8$, we have
$\jump[2]{E_8}{M(A_4)}=24$.

Apart from trivialities, the remaining constants to be calculated are
\begin{eqnarray*}
    &&\jump[i]{\cyclic{2}}{A_4},\;\;
        \jump[i]{S_3}{C},\;\;
        \jump[i]{A_4}{C},\;\;
        \jump[i]{E_4^-}{C},\\
    &&\jump[i]{E_4^-}{M(A_4)},\;\;
        \jump[i]{E_4^+}{M(A_4)},\;\;
        \jump[i]{E_4^+}{C},\;\;
        \jump[i]{E_8}{C},
\end{eqnarray*}
for $i=2$, $3$. As some invention is needed here, we show how to obtain all of
them.

We begin with $\jump[i]{S_3}{C}$. Let $G$ be a copy of $S_3$ in $C$. For any
element $x\not\in G$, we must have $\spn{G,\,x}\cong M(S_3)$, $M(A_4)$, or $C$.
Therefore, $\jump[i]{S_3}{C} = (i-1)\cdot\glbidx{\cyclic{i}}{C} -
\jump[i]{S_3}{M(S_3)} - \jump[i]{S_3}{M(A_4)} -
(i-1)\cdot\glbidx{\cyclic{i}}{S_3}$, for $i=2$, $3$. Consequently,
\begin{displaymath}
    \jump[2]{S_3}{C}=63-6-36-3=18,\;\;\;
    \jump[3]{S_3}{C}=56-0-18-2=36.
\end{displaymath}
Similarly,
\begin{displaymath}
    \begin{array}{rclrcl}
    \jump[2]{\cyclic{2}}{A_4}&=&0,&\jump[3]{\cyclic{2}}{A_4}&=&24,\\
    \jump[2]{A_4}{C}&=&48,&\jump[3]{A_4}{C}&=&48,\\
    \jump[2]{E_8}{C}&=&32,&\jump[3]{E_8}{C}&=&32.
    \end{array}
\end{displaymath}

A more detailed analysis of the subloop lattice of $C$ allows us to calculate
the remaining eight constants.

\begin{figure}
    \centering
    \input{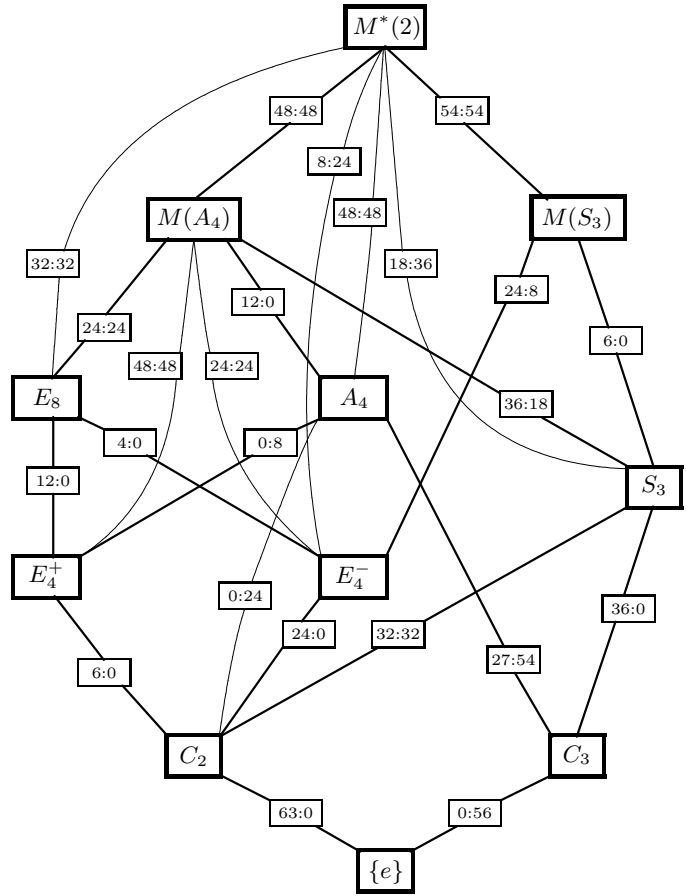}
    \caption[]{The constants $\jump[i]{A}{B}$ for $\paige{2}$. If $A$ is
    maximal in a copy of $B$, then $A$ and $B$ are connected by a thick,
    straight line; else by a thin, curved line. The constants $\jump[2]{A}{B}$,
    $\jump[3]{A}{B}$ are located in a box on the edge connecting $A$ and $B$,
    separated by colon.}
    \label{Fg:Random}
\end{figure}

\begin{lemma}
Let $G\in\orbit[]^-$, and let $M_1$, $M_2$, $M_3$ be the three copies of
$M(A_4)$ containing $G$. Then $M_i\cap M_j$ contains no element of order $3$,
for $i\ne j$, and $M_1\cap M_2\cap M_3$ is the unique copy of $E_8$ containing
$G$. In particular,
\begin{displaymath}
    \jump[3]{E_4^-}{M(A_4)}=24,\;\;\jump[3]{E_4^-}{C}=24,\;\;
    \jump[2]{E_4^-}{M(A_4)}=24,\;\;\jump[2]{E_4^-}{C}=8.
\end{displaymath}
\begin{proof}
Assume there is $x\in M_i\cap M_j$, $\order{x}=3$, for some $i\ne j$. Then
$M_i=\spn{G,\,x}=M_j$, because $\isoidx{E_4^{-}}{A_4}{C}=0$, a contradiction.
Thus $M_1\cup M_2\cup M_3$ contains $3\cdot 8=24$ elements $x$ of order $3$
such that $\spn{G,\,x}\in\orbit[M(A_4)]$.

Let $H$ be the unique copy of $E_8$ containing $G$. We must have $H=M_1\cap
M_2\cap M_3$, since $\isoidx{E_8}{M(A_4)}{C}=3$. Therefore $M_1\cup M_2\cup
M_3$ contains $3\cdot(12-4)=24$ involutions $x$ such that
$\spn{G,\,x}\in\orbit[M(A_4)]$. The constants $\jump[i]{E_4^-}{C}$ are then
easy to calculate with help of Figure \ref{Fg:Random}.
\end{proof}
\end{lemma}

It is conceivable that there is $G\in\orbit[]^+$ and $x\in C$ such that
$\spn{G,\,x}=C$. It is not so, though.

\begin{lemma}\label{Lm:MoreConstants}
In $C$, we have
\begin{displaymath}
    \jump[3]{E_4^+}{M(A_4)}=48,\;\;\jump[2]{E_4^+}{M(A_4)}=48,\;\;
    \jump[3]{E_4^+}{C}=0,\;\;\jump[2]{E_4^+}{C}=0.
\end{displaymath}
\begin{proof}
Pick $G\in\orbit[]^+$, and let $M_1$, $\dots$, $M_7$ be the seven copies of
$M(A_4)$ containing $G$. We claim that $(M_i\cap M_j)^2=\neutral$, for $i\ne
j$. Assume it is not true, and let $x$ be an element of order $3$ contained in
$M_i\cap M_j$. Then $A_4\cong\spn{G,\,x}\le M_i\cap M_j$ shows that
$\glbidx{A_4}{M(A_4)}\ge 2$, a contradiction. Thus $\bigcup_{i=1}^7 M_i$
contains all $8\cdot 7=56$ elements of order $3$. In particular, we have
$\spn{G,\,x}\ne C$ for any element $x$ of order $3$. This translates into
\begin{displaymath}
    \jump[3]{E_4^+}{C}=0,\;\;
    \jump[3]{E_4^+}{M(A_4)}=56-\jump[3]{E_4^+}{A_4}=48.
\end{displaymath}

We proceed carefully to show that $\jump[2]{E_4^+}{C}=0$. The group $G$ is
contained in a single copy $A$ of $A_4$, that is in turn contained in a single
copy $M_1$ of $M(A_4)$. Let $H_1$, $H_2$, $H_3\le M_1$ be the three copies of
$E_8$ containing $G$ (see the proof of Lemma 13.1 \cite{VojtechovskyHasse}).
Observe that $H_1\cup H_2\cup H_3=G\cup Au$, where $Au$ is the second coset of
$A$ in $M_1$. Pick $M_i$, $M_j$, with $2\le i<j\le 7$. We want to show that
$M_i\cap M_j\subseteq M_1$. Thanks to the first part of this Lemma, we know
that $M_i\cap M_j\cong E_4$ or $E_8$. When $M_i\cap M_j\cong E_4$ then,
trivially, $M_i\cap M_j=G\le M_1$. When $M_i\cap M_j\cong E_8$ then $M_i\cap
M_j=H_k$ for some $k\in\{1,\,2,\,3\}$, else $\isoidx{G}{E_8}{C}\ge 4$, a
contradiction.

Consequently, $\bigcup_{i=1}^7 M_i$ contains at least $15+6\cdot 8=63$
involutions; $15$ in $M_1$, and additional $8$ in each $M_i$, $i>1$. In
particular, $\spn{G,\,x}\ne C$ for every involution $x$. We get
\begin{displaymath}
    \jump[2]{E_4^+}{C}=0,\;\;
    \jump[2]{E_4^+}{M(A_4)}
    =60-\jump[2]{E_4^+}{E_8}=48.
\end{displaymath}
This finishes the proof.
\end{proof}
\end{lemma}

\setlength{\unitlength}{0.85mm}

All constants $\jump[i]{A}{B}$ have now been calculated. They are collected in
Figure \ref{Fg:Random}.

\subsection{Random generators of arbitrary orders}

\begin{figure}
    \centering
    \input{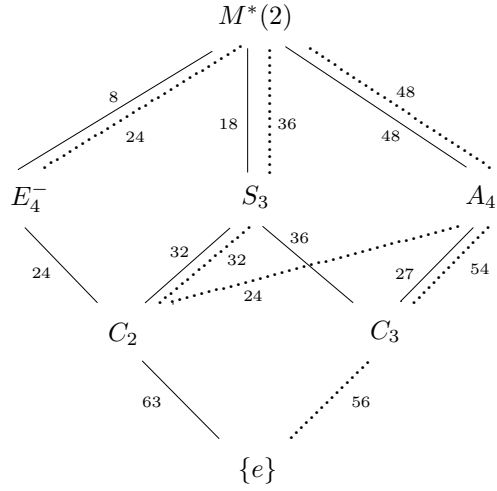}
    \caption[]{Sequences of subloops in $M^*(2)$}
    \label{Fg:AllSequences}
\end{figure}

\noindent According to Figure \ref{Fg:Random}, there are only five
orbit-nonequivalent ways to get from $\neutral$ to $C$ in $3$ steps. Namely,
\begin{eqnarray*}
    \mathbf A_0&=&\{\{\neutral\},\,\cyclic{2},\,A_4,\,C\},\\
    \mathbf A_1&=&\{\{\neutral\},\,\cyclic{2},\,E_4^-,\,C\},\\
    \mathbf A_2&=&\{\{\neutral\},\,\cyclic{2},\,S_3,\,C\},\\
    \mathbf A_3&=&\{\{\neutral\},\,\cyclic{3},\,S_3,\,C\},\\
    \mathbf A_4&=&\{\{\neutral\},\,\cyclic{3},\,A_4,\,C\}.
\end{eqnarray*}
These sequences and the related constants $\jump[i]{A}{B}$ are visualized in
Figure \ref{Fg:AllSequences}. Full lines correspond to involutions $(i=2)$,
dotted lines to elements of order $3$ $(i=3)$.

\setlength{\unitlength}{1.0mm}

\begin{proposition}
Let $C=\paige{2}$. Then the probability that $3$ randomly chosen elements of
$C$ generate $C$ is $\prob{3}{C} = 955,584\cdot 120^{-3} = 0.553$.
\begin{proof}
By $(\ref{Eq:Gen})$,
\begin{displaymath}
    |\gen{3}{C}|=\sum_{i=0}^4|\orbclass[]{\mbf{A_i}}|.
\end{displaymath}
By our previous calculation summarized in Figure \ref{Fg:AllSequences},
$|\orbclass[]{\mbf{A_0}}|=63\cdot 24\cdot (48+48)$,
$|\orbclass[]{\mbf{A_1}}|=63\cdot24\cdot(8+24)$,
$|\orbclass[]{\mbf{A_2}}|=63\cdot(32+32)\cdot(18+36)$,
$|\orbclass[]{\mbf{A_3}}|=56\cdot 36\cdot(18+36)$, and
$|\orbclass[]{\mbf{A_4}}|=56\cdot(54+27)\cdot(48+48)$. Thus
$|\gen{3}{C}|=955,584$. We are done by $(\ref{Eq:Prob})$.
\end{proof}
\end{proposition}

\subsection{Random generators of given orders}

\begin{figure}
    \centering
    \input{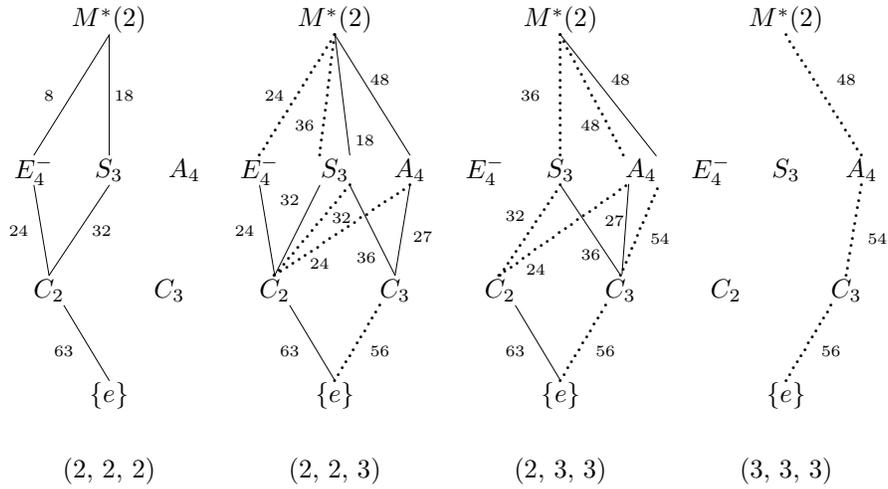}
    \caption[]{The shortest sequences of subloops in $M^*(2)$}
    \label{Fg:SeparateSequences}
\end{figure}

\noindent The only possible types of orders for three generators in $C$ are
$(2,\,2,\,2)$, $(2,\,2,\,3)$, $(2,\,3,\,3)$, and $(3,\,3,\,3)$. The sequences
of subloops corresponding to each of these types are depicted in Figure
\ref{Fg:SeparateSequences}. We must be careful, though, since not all
combinations of lines in Figure \ref{Fg:SeparateSequences} correspond to
sequences with correct types of orders. The possible continuations are
emphasized in Figure \ref{Fg:SeparateSequences}.

\setlength{\unitlength}{1mm}

\begin{proposition}
Let $\mathbf s=(s_1,\,s_2,\,s_3)$ be a $3$-tuple of integers, $s_1\le s_2\le
s_3$, $C=\paige{2}$, and let $\prob{\mathbf s}{C}$ be the probability that $3$
randomly chosen elements $a_1$, $a_2$, $a_3$ of $C$ generate $C$ and
$(\order{a_1},\,\order{a_2},\,\order{a_3})$ is of type $\mathbf s$. Then
$\prob{(2,\,2,\,2)}{C} = 48,384\cdot 120^{-3} = 0.028$, $\prob{(2,\,2,\,3)}{C}
= 326,592\cdot 120^{-3} = 0.189$, $\prob{(2,\,3,\,3)}{C}=435,456\cdot 120^{-3}
= 0.252$, and $\prob{(3,\,3,\,3)}{C} = 145,152\cdot 120^{-3} = 0.084$.
\begin{proof}
Use Figure \ref{Fg:SeparateSequences} and $(\ref{Eq:Tprob})$.
\end{proof}
\end{proposition}

\bibliographystyle{plain}

\begin{thebibliography}{9}

\bibitem[Kantor and Seress(2001)]{GroupsAndComputationIII}
Groups and computation III, proceedings of the 3rd International Conference
held at The Ohio State University, Columbus, OH, June $15$--$19$, $1999$.
Edited by William M.~Kantor and \'Akos Seress. Ohio State University
Mathematical Research Institute Publications, \textbf{8}, Walter de Gruyter,
Berlin, 2001.

\bibitem[Vojt\v echovsk\'y(2001)]{VojtechovskyHasse} P.~Vojt\v echovsk\'y,
\emph{Investigation of subalgebra lattices by means of Hasse constants}, to
appear in Algebra Universalis.

\end{thebibliography}

\end{document}